\documentclass[11pt]{article}

\usepackage{latexsym,amssymb}

\makeatletter \@addtoreset{equation}{section}

\renewcommand{\Re}{{\rm Re}\,}

\newtheorem{thm}{Theorem}[section]
\newtheorem{prop}[thm]{Proposition}
\newtheorem{lem}[thm]{Lemma}
\newtheorem{cor}[thm]{Corollary}

\def\enu#1{\newline\makebox[5mm][l]{\rm(#1)}}
\def\benu#1#2{\newline\makebox[#1 mm][l]{\rm(#2)}}

\def\bp{\noindent{\it Proof.}\ }
\def\bpp#1{\smallskip\noindent{\it Proof of #1.}\ }
\def\ep{\nopagebreak\newline\mbox{\ }\hfill\rule{2mm}{2mm}}
\def\epp{\nopagebreak\mbox{\ }\hfill\rule{2mm}{2mm}}

\def\Ch{{\rm Ch}}
\def\End{{\rm End}}

\def\ind{{\rm Ind}_F}

\def\indq{q{\rm-Ind}_F}

\def\Ker{{\rm Ker}}

\def\span{{\rm span}}

\def\Tr{{\rm Tr}}

\def\7#1{{\mathbb #1}}

\def\A{{\cal A}}
\def\B{{\cal B}}
\def\C{{\cal C}}

\def\H{{\cal H}}

\def\eps{{\varepsilon}}

\def\rad{{\blacktriangleleft}}
\def\r{{\triangleleft}}

\def\pr{{\partial}}

\def\2{{\frac{1}{2}}}

\def\>{{\rangle}}
\def\<{{\langle}}

\def\D{{\Delta}}

\def\Ir0{{\rm Irr}(A_0,\Delta_0)}
\def\sl2{{U_q({\mathfrak su}_2)}}

\begin{document}

\title{\bf A Local Index Formula for the Quantum Sphere}

\author{Sergey Neshveyev$^1$
\& Lars Tuset$^2$}

\footnotetext[1]{Partially supported by the Norwegian Research
Council.}

\footnotetext[2]{Supported by the SUP-program of the Norwegian
Research Council.}

\date{}

\maketitle

\begin{abstract}
For the Dirac operator $D$ on the standard quantum sphere we
obtain an asymptotic expansion of the $SU_q(2)$-equivariant entire
cyclic cocycle corresponding to $\eps^\2 D$ when evaluated on the
element $k^2\in\sl2$. The constant term of this expansion is a
twisted cyclic cocycle which up to a scalar coincides with the
volume form and computes the quantum as well as the classical
Fredholm indices.
\end{abstract}

\section*{Introduction}

The geometry of $q$-deformed spaces and the theory of
non-commutative geometry of Connes have developed considerably
over the last twenty years. Crucial in establishing connections
between these two areas is the construction of Dirac operators which
give rise to reasonable differential calculi. Other aspects of
the theory include computations of cyclic cohomology and analysis
of the index theorem. Work in this direction has progressed
furthest for spaces like $SU_q(2)$ and the quantum spheres, see e.g.
\cite{CP,DS,SW,MNW1,MNW2,Co2}.

In this paper we consider the index formula for the homogeneous
sphere $S^2_q$ of Podle\'s. The Dirac operator for this space has
very different properties compared to the classical one. In
particular, the associated $\zeta$-function has infinitely many poles
on vertical lines and
the traces $\Tr(e^{-tD^2})$ of the heat operators tend to infinity
slower than $t^{-p}$ for any $p>0$. More importantly, the spectral
triple defined by the Dirac operator does not satisfy the
regularity assumption, a condition which is often overshadowed by
other assumptions such as boundedness of commutators, but which is
crucial for a detailed analysis. This can mean that the Dirac
operator is, in fact, not the right one, and one can try to
construct another operator  with the same Fredholm module.
However, in this process one will most likely lose the
$SU_q(2)$-equivariance~\cite{DS}. What one gains is not so clear:
with the absence of Getzler's symbol calculus and of a
Wodzicki-type geometric description of the Dixmier trace, the
computation of the indices in terms of the Dixmier trace remains
for the moment a non-trivial problem. So we will stick to the
standard Dirac operator. The associated spectral triple is
$SU_q(2)$-equivariant and thus, via the JLO-cocycle, defines the
Chern character in the equivariant entire cyclic cohomology of
$S^2_q$. We evaluate this cocycle on $k^2\in\sl2$. The resulting
object  is an entire twisted cyclic cocycle, and its pairing with
the equivariant $K$-theory computes the quantum Fredholm index,
that is, the difference of quantum dimensions of the kernel and
the cokernel~\cite{NT}. The possibility of detecting cocycles by
evaluating them on $k^2$ was pointed out by Connes in~\cite{Co2}.
As was promised in~\cite{NT} such evaluations are much easier to
compute. The philosophical reason for this is that the quantum
dimension is intrinsically associated with the tensor category of
finite dimensional corepresentations of a compact quantum group.
On the technical side, the evaluation on $k^2$ gives an immediate
connection to the Haar state, as was already remarked
in~\cite{G,SW}, and thus serves as a replacement for the trace
theorem by Connes on equality of the Dixmier trace and the
Wodzicki residue. Note also that for $q$-deformations the
classical index can be recovered from the quantum one: if one can
prove that the Fredholm operators depend continuously on $q$ and
has a polynomial formula in $q$ and $q^{-1}$ for the quantum
index, the classical index is obtained by simply setting $q=1$.

The evaluation on $k^2$ does not, however, solve all problems: the
spectral triple is still non-regular, and while the quantum traces
of the heat operators are better controlled than the classical
ones, they possess a strange oscillating behavior near zero. The
non-regularity of the spectral triple can be illuminated by saying
that though the principal symbol of $|D|$ is scalar-valued, the
operator $|D|^zT|D|^{-z}$ does not necessarily have the same
symbol as $T$. The symbol is, however, computable, and this is the
key property allowing to apply the ideas of the proof of the local
index formula of Connes-Moscovici~\cite{CM3} to our situation. It
seems that the development of a full-scaled pseudo-differential
calculus for the Dirac operators on $q$-deformed spaces is the
main prerequisite for the analysis of more general examples, such
as the ones in~\cite{Kr}.

The Dirac operator on the quantum sphere is defined using the
standard differential calculus, and in the course of our analysis
we use very little beyond basic properties of this calculus and
some information on the spectrum of the Dirac operator. In this
respect the fact that Connes' non-commutative geometry can be
applied successfully to the study of a $q$-deformed space is more
important than the final result for the quantum sphere.

{\bf Acknowledgement} The preparation of this paper was finished
during the authors' stay at Institute Mittag-Leffler in September
2003. They would like to express their gratitude to the staff at
the institute and to the organizers of the year in "Noncommutative
Geometry".

\bigskip

\section{Cyclic Cohomology} \label{1}

In this section we formulate some results of non-commutative
geometry in the Hopf algebra equivariant setting. Fortunately, as
the Hopf algebra equivariant cyclic cohomology theory is now
available \cite{AK1,NT}, a simple argument involving crossed
products allows to transfer such results from the non-equivariant
case without having to prove each one of them from scratch.

We use the same notation and make the same assumptions as in 
\cite{NT}. So let $(\H,\D)$
be a Hopf algebra over $\7C$ with invertible antipode $S$ and
counit $\eps$, and adapt the Sweedler notation
$\D(\omega)=\omega_{(0)}\otimes\omega_{(1)}$. The results below
are also true for the algebra of finitely supported functions on a
discrete quantum group, which is of main interest for us. Let $\B$
be a unital right $\H$-module algebra with right action of
$(\H,\D)$ denoted by~$\r$. Consider the space $C^n_\H(\B)$ of
$\H$-invariant $n$-cochains, so $C^n_\H(\B)$ consists of linear
functionals on $\H\otimes\B^{\otimes(n+1)}$ such that
$$
f(S^{-1}(\omega_{(0)})\eta\omega_{(1)};b_0\r\omega_{(2)},\ldots,
b_n\r \omega_{(n+2)})
=\eps(\omega)f(\eta;b_0,\ldots,b_n)
$$
for any $\omega,\eta\in\H$ and $b_0,\ldots,b_n\in\B$. The
coboundary operator $b_n\colon C^{n-1}_\H(\B)\to C^n_\H(\B)$ and
the cyclic operator $\lambda_n\colon C^n_\H(\B)\to C^n_\H(\B)$ are
given by
$$
(b_nf)(\omega;b_0,\ldots,b_n)=\sum^{n-1}_{i=0}(-1)^i
f(\omega;b_0,\ldots,b_{i-1},b_ib_{i+1},b_{i+2},\ldots,b_n)
$$
$$
\hspace{3cm}+(-1)^nf(\omega_{(0)};b_n\r
\omega_{(1)}b_0,b_1,\ldots,b_{n-1})
$$
and
$$
(\lambda_nf)(\omega;b_0,\ldots,b_n)
=(-1)^nf(\omega_{(0)};b_n\r\omega_{(1)},b_0,\ldots,b_{n-1}).
$$
The subcomplex $(\Ker(\iota-\lambda),b)$ of $(C^\bullet_\H(\B),b)$
is denoted by  $(C^\bullet_{\H,\lambda}(\B),b)$ and its cohomology
is denoted by $HC^\bullet_\H(\B)$. There is a pairing
$\<\cdot,\cdot\>\colon HC^{2n}_\H(\B)\times K^\H_0(\B)\to R(\H)$
of the even cyclic cohomology with equivariant $K$-theory, where
$R(\H)$ is the space of $\H$-invariant linear functionals on $\H$
with $\H$ acting on itself by
$$
\eta\r\omega=S^{-1}(\omega_{(0)})\eta\omega_{(1)}.
$$
If $p\in\B$ is an $\H$-invariant idempotent and $f\in
C^{2n}_{\H,\lambda}(\B)$ is a cyclic cocycle, then
\begin{equation} \label{e1.1}
\<[f],[p]\>(\omega)=\frac{1}{n!}f(\omega;p,\ldots,p).
\end{equation}
More generally, if $X$ is a finite dimensional (as a vector space)
right $\H$-module, then $\End(X)\otimes\B$ is a right $\H$-module
algebra with right action $\rad$ given by
$$
(T\otimes b)\rad\omega
=\pi_X(\omega_{(0)})T\pi_XS^{-1}(\omega_{(2)})\otimes b\r\omega_{(1)},
$$
where $\pi_X\colon\H\to\End(X)$ is the anti-homomorphism defining
the $\H$-module structure on $X$. Then we have a map
$\Psi^n_X\colon C^n_\H(\B)\to C^n_\H(\End(X)\otimes\B)$ given by
$$
(\Psi^n_Xf)(\omega;T_0\otimes b_0,\ldots,T_n\otimes b_n)
=f(\omega_{(0)};b_0,\ldots,b_n)
\Tr(\pi_XS^{-1}(\omega_{(1)})T_0\ldots T_n).
$$
Now suppose $p\in\End(X)\otimes\B$ is an $\H$-invariant
idempotent, so it defines an element of $K^\H_0(\B)$. By
definition we have
$$
\<[f],[p]\>=\<[\Psi^{2n}_Xf],[p]\>
$$
for a cyclic cocycle $f\in
C^{2n}_{\H,\lambda}(\B)$.

\smallskip

We will also need the $(b,B)$-bicomplex description of the
periodic cyclic cohomology. So consider the operator $B=NB_0$,
where $N_n\colon C^n_\H(\B)\to C^n_\H(\B)$,
$N_n=\sum^n_{i=0}\lambda_n^i$, $B^n_0\colon C^{n+1}_\H(\B)\to
C^n_\H(\B)$, $B_0^n=(-1)^ns^n_n(\iota-\lambda_n)$, so
$$
(B^n_0f)(\omega;b_0,\ldots,b_n)=f(\omega;1,b_0,\ldots,b_n)-
(-1)^{n+1}f(\omega;b_0,\ldots,b_n,1).
$$
Set $\C^0_\H(\B)=\oplus^\infty_{n=0}C^{2n}_\H(\B)$ and
$\C^1_\H(\B)=\oplus^\infty_{n=0}C^{2n+1}_\H(\B)$. Then we have a
well-defined complex
$$
\C^1_\H(\B)\stackrel{b+B}\longrightarrow\C^0_\H(\B)
\stackrel{b+B}\longrightarrow\C^1_\H(\B)
\stackrel{b+B}\longrightarrow\C^0_\H(\B).
$$
The cohomology of this small complex is denoted by $HP^n_\H(\B)$,
$n=0,1$. Then the formula
\begin{equation} \label{e1.2}
\<[(f_{2n})_n],[p]\>(\omega)
=\sum^\infty_{n=0}(-1)^n\frac{(2n)!}{n!}
(\Psi^{2n}_Xf_{2n})(\omega;p-\2,p,\ldots,p)
\end{equation}
defines a pairing $HP^0_\H(\B)\times K^\H_0(\B)\to R(\H)$ (this
can be deduced from the non-equivariant case in the same way as in
\cite{NT}). The map
\begin{equation} \label{e1.3}
C^{2n}_\H(\B)\to\C^0_\H(\B),\ \ f\mapsto\frac{(-1)^n}{(2n)!}f,
\end{equation}
induces a homomorphism $HC^{2n}_\H(\B)\to HP^0_\H(\B)$ compatible
with pairings (\ref{e1.1}) and (\ref{e1.2}) and which respects the
periodicity operator in the sense that the images of $f$ and $Sf$
in $HP^0$ are equal.

\smallskip

Consider the subspaces
$CE^k_\H(\B)\subset\prod^\infty_{n=0}C^{2n+k}_\H(\B)$, $k=0,1$,
consisting of all cochains $(f_{2n+k})_n$ such that for any
$\omega\in\H$ and any finite subset $F\subset\B$ there exists
$C>0$ such that
$$
|f_{2n+k}(\omega_{(0)};b_0\r\omega_{(1)},\ldots,
b_{j-1}\r\omega_{(j)},b_j,\ldots,b_{2n+k})|
\le\frac{C}{n!}
$$
for any $n\ge0$, $0\le j\le 2n+k+1$ and $b_i\in F$. Denote by
$HE^k_\H(\B)$ the cohomology of the complex
$$
CE^1_\H(\B)\stackrel{b+B}\longrightarrow CE^0_\H(\B)
\stackrel{b+B}\longrightarrow CE^1_\H(\B)
\stackrel{b+B}\longrightarrow CE^0_\H(\B).
$$
The formula (\ref{e1.2}) still defines a pairing of $HE^0_\H(\B)$
with $K^\H_0(\B)$.

\medskip

Suppose we are given an equivariant even spectral triple. By
this we mean a Hilbert space $H$ with grading operator $\gamma$,
an odd selfadjoint operator $D$ with compact resolvent, and even
representations of $\B$ and $\H$ on $H$ such that $D$ commutes
with $\H$, $[D,b]$ is bounded for any $b\in\B$ and
\begin{equation} \label{e1.4}
b\omega=\omega_{(0)}b\r\omega_{(1)}\ \ \forall b\in\B,\ \forall
\omega\in\H.
\end{equation}
If the spectral triple is  $\theta$-summable in the sense that
$e^{-tD^2}$ is trace-class for all $t>0$, set
$$
\Ch^{2n}(D)(\omega;b_0,\ldots,b_{2n})=\int_{\D_{2n}}dt \Tr(\gamma
\omega
b_0e^{-t_0D^2}[D,b_1]e^{-t_1D^2}\ldots[D,b_{2n}]e^{-t_{2n}D^2}),
$$
where $\int_{\D_{2n}}dt$ means integration over the simplex
$\D_{2n}=\{(t_0,\ldots,t_{2n})\,|\,t_i\ge0,\sum_it_i=1\}$. Then
$(\Ch^{2n}(D))_n$ is a cocycle in $CE^0_\H(\B)$ having all the
usual properties of a JLO-cocycle~\cite{Co1,CM2,GBVF}:
\vspace{2mm}\enu{i} $(\Ch^{2n}(tD))_n$ and $(\Ch^{2n}(sD))_n$ are
cohomologous for any $s,t>0$; \vspace{2mm}\enu{ii} the pairing of
$K^\H_0(\B)$ with $(\Ch^{2n}(tD))_n$ computes the index map
defined by $D$ (so e.g. for an $\H$-invariant idempotent $p\in B$
we get $\<[(\Ch^{2n}(D))_n],[p]\>(\omega)=\Tr(\omega|_{\Ker\,
p_-Dp_+})-\Tr(\omega|_{\Ker\, p_+Dp_-})$);
\vspace{2mm}\benu{6}{iii} if the spectral triple is $p$-summable
in the sense that $|D|^{-p}$ is a trace-class operator,
$F=D|D|^{-1}$, and $\tau^{2m}_F$ denotes the Chern character of
the Fredholm module $(H,F,\gamma)$ in the equivariant cyclic
cohomology for some fixed $m$, $2m\ge p$, so
$$
\tau^{2m}_F(\omega;b_0,\ldots,b_{2m})=\frac{(-1)^m}{2}m!
\Tr(\gamma\omega F[F,b_0]\ldots[F,b_{2m}]),
$$
then $(\Ch^{2n}(D))_n$ is cohomologous to the image of
$\tau^{2m}_F$ under the map $C^{2m}_\H(\B)\to CE^0_\H(\B)$ given
by (\ref{e1.3}).

As in the non-equivariant case, all these properties are
consequences of the homotopy invariance of the cohomology class of
$(\Ch^{2n}(D))_n$ meaning that if we have an $\H$-invariant
homotopy $D_t$ such that e.g. $\dot D_t$ is bounded, then
$(\Ch^{2n}(D_0))_n$ and $(\Ch^{2n}(D_1))_n$ are cohomologous
(concerning (ii) see also~\cite{KL}). This in turn can be deduced
from the non-equivariant case as follows, see the discussion at
the end of \cite{Co2}. Consider the crossed product algebra
$\B\rtimes\H$, i.e. the vector space $\B\otimes\H$ with product
$$
(b\otimes\omega)(c\otimes\eta)=bc\r
S^{-1}(\omega_{(1)})\otimes\omega_{(0)}\eta.
$$
In the sequel we write $b\omega$ instead of $b\otimes\omega$. Then
we have a map $\Phi^n\colon C^n_\H(\B)\to C^n(\B\rtimes\H)$ given
by
$$
(\Phi^nf)(b_0\omega^0,\ldots,b_n\omega^n)
=f(\omega^0_{(0)}\ldots\omega^n_{(0)};
b_0\r(\omega^0_{(1)}\ldots\omega^n_{(1)}),
b_1\r(\omega^1_{(2)}\ldots\omega^n_{(2)}),
\ldots,b_n\r\omega^n_{(n+1)}).
$$
The definition of $\Phi^n$ is motivated by the following easily
proved lemma.

\begin{lem} \label{1.1}
Consider a covariant representation of $\B$ and $\H$ on $H$
(that is, identity (\ref{e1.4}) holds). Suppose $c_0,\ldots,c_n\in
B(H)$ commute with $\H$. Then for any $b_0,\ldots,b_n\in\B$ and
$\omega^0,\ldots,\omega^n\in\H$ we have
$$
c_0b_0\omega^0\ldots c_nb_n\omega^n
=\omega^0_{(0)}\ldots\omega^n_{(0)}c_0
b_0\r(\omega^0_{(1)}\ldots\omega^n_{(1)})c_1
b_1\r(\omega^1_{(2)}\ldots\omega^n_{(2)}) \ldots
c_nb_n\r\omega^n_{(n+1)}.
$$
\epp
\end{lem}
It is not difficult to check that the maps $\Phi^n$ constitute a
morphism of cocyclic objects, so they induce maps for the various
cyclic cohomology theories. Note also that if $\H$ is unital, or
$\H$ has an approximate unit in an appropriate sense, then
$\Phi^n$ is injective. So to prove a property for a cochain in an
equivariant theory it is enough to establish the analogous
property for the image of the cochain under $\Phi$.
%One has only
%to be sure that constructions appearing in the proof are
%equivariant.
Now note that any equivariant spectral triple can be considered as
a spectral triple for $\B\rtimes\H$. It follows immediately from
Lemma~\ref{1.1} that $(\Phi^{2n}(\Ch^{2n}(D))_n$ is the
JLO-cocycle for this spectral triple. This allows to deduce the
properties of $(\Ch^{2n}(D))_n$ from the non-equivariant case.

\medskip

From now onwards we assume that $(\H,\D)$ is the algebra of
finitely supported functions on a discrete quantum group $(\hat
A,\hat\D)$, and as in \cite{NT} we write $(\hat\A,\hat\D)$ instead
of $(\H,\D)$. Denote by $\rho\in M(\hat\A)$ the Woronowicz
character $f_{-1}$. We will be interested in evaluating
equivariant cocycles on $\omega=\rho$. The following known lemma
explains why such evaluations are easier to deal with.

\begin{lem} \label{1.2}
Let $U\in M(A\otimes K(H))$ be a unitary corepresentation of the
dual compact group $(A,\D)$, and $\alpha_U\colon B(H)\to M\otimes
B(H)$ the coaction of the von Neumann closure $(M,\D)$ of $(A,\D)$
on~$B(H)$,
$$
\alpha_U(x)=U^*(1\otimes x)U.
$$
Consider also the corresponding representation of $\hat A$ on $H$,
so $\omega\xi=(\omega\otimes\iota)(U)\xi$. Then the map
$a\mapsto\Tr(\cdot a\rho)$ defines a one-to-one correspondence
between positive elements $a\in B(H)^{\alpha_U}$ such that
$\Tr(a\rho)=1$ and normal $\alpha_U$-invariant states on $B(H)$.

In particular, if $B\subset B(H)$ is an $\alpha_U$-invariant
C$^*$-subalgebra with a unique $\alpha_U$-invariant
state~$\varphi$, then
$$
\Tr(ba\rho)=\varphi(b)\Tr(a\rho)
$$
for any $b\in B$ and any $a\in B(H)^{\alpha_U}_+$ with
$\Tr(a\rho)<\infty$. \epp
\end{lem}

Note also that the fixed point algebra $B(H)^{\alpha_U}$ is
precisely the commutant of $\hat A$.

Since in general $\rho\notin\hat\A$, one needs some care in
dealing with this element. From now on we assume as in
\cite{NT} that we have a left coaction $\alpha$ of $(A,\D)$ on a
C$^*$-algebra $B$, and $\r$ is the corresponding action of
$\hat\A$,
$$
b\r\omega=(\omega\otimes\iota)\alpha(b),
$$
while $\B=B\r\hat\A$. Then $\sigma_z(b)=b\r\rho^z$ is a
well-defined element of $\B$ for any $b\in\B$ and $z\in\7C$. The
automorphism $\sigma=\sigma_1$ will be called the twist. Now if we
have an equivariant $\theta$-summable spectral triple, for the
expression $\Ch^{2n}(D)(\rho;b_0,\ldots,b_{2n})$ to make sense it
is enough to require that
$$
\Tr(\rho e^{-sD^2})<\infty\ \ \forall s>0.
$$
Indeed, since $z\mapsto\rho^{-z}[D,b]\rho^z=[D,b\r\rho^z]\in B(H)$
is an analytic function, the H\"{o}lder inequality and the
identity

\medskip\noindent$\displaystyle
e^{-t_0D^2}[D,b_1]e^{-t_1D^2}\ldots[D,b_{2n}]e^{-t_{2n}D^2}\rho$
\begin{equation} \label{e1.10}
=e^{-t_0D^2}\rho^{t_0}[D,b_1\r\rho^{t_0}]e^{-t_1D^2}\rho^{t_1}
[D,b_2\r\rho^{t_0+t_1}]e^{-t_2D^2}\rho^{t_2} \ldots
[D,b_{2n}\r\rho^{1-t_{2n}}]e^{-t_{2n}D^2}\rho^{t_{2n}}
\end{equation}
show that $\Ch^{2n}(D)(\rho;b_0,\ldots,b_{2n})$ is well-defined.
In fact, these elements form a cocycle in an appropriate cyclic
cohomology theory, namely, in the entire twisted cyclic cohomology
$HE^0_\sigma(\B)$, cf \cite{KMT,G}. This is immediate as such
twisted theories are obtained by setting $\omega=\rho$ in the
formulas above. Thus the space of $n$-cochains $C^n_\sigma(\B)$
becomes the space of linear functionals $f$ on $\B^{\otimes(n+1)}$
such that $f(\sigma(b_0),\ldots,\sigma(b_n)) =f(b_0,\ldots,b_n)$,
the cyclic operator is given by
$(\lambda_nf)(b_0,\ldots,b_n)=(-1)^nf(\sigma(b_n),b_0,\ldots,b_{n-1})$,
and so on. In particular, we still have pairings of $K^\H_0(\B)$
with $HC^{2n}_\sigma(\B)$, $HP^0_\sigma(\B)$ and
$HE^0_\sigma(\B)$.

\bigskip

\section{Differential Calculus on the Quantum Sphere} \label{2}

From now on $(A,\D)$ will denote the compact quantum group $SU_q
(2)$ of Woronowicz \cite{Wor}, $q\in(0,1)$. So $A$ is the
universal unital C$^*$-algebra with generators $\alpha$ and
$\gamma$ satisfying the relations
$$
\alpha^*\alpha+\gamma^*\gamma =1, \ \ \alpha\alpha^*
+q^2\gamma^*\gamma =1,\ \ \gamma^*\gamma =\gamma\gamma^* , \ \
\alpha\gamma =q\gamma\alpha,\ \ \alpha\gamma^*=q\gamma^*\alpha .
$$
The comultiplication $\D$ is determined by the formulas
$$
\D (\alpha )=\alpha\otimes\alpha -q\gamma^*\otimes\gamma,\ \
\D (\gamma )=\gamma\otimes\alpha +\alpha^*\otimes\gamma .
$$

Consider the quantized universal enveloping algebra $\sl2\subset
M(\hat\A)$. Recall that it is the universal unital $*$-algebra
generated by elements $e,f,k,k^{-1}$ satisfying the relations
$$
kk^{-1}=k^{-1}k=1,\ \ ke=qek,\ \ kf=q^{-1}fk,\ \
ef-fe=\frac{k^2-k^{-2}}{q-q^{-1}},
$$
$$
k^*=k,\ \ e^* =f.
$$
The algebra $\sl2$ is a Hopf subalgebra of $M(\hat\A)$ with
$$
\hat\D (k)=k\otimes k,\ \
\hat\D (e)=e\otimes k^{-1}+k\otimes e,\ \
\hat\D (f)=f\otimes k^{-1}+k\otimes f,
$$
$$
\hat S(k)=k^{-1},\ \ \hat S(e)=-q^{-1}e,\ \ \hat S(f)=-qf,
$$
$$
\hat\varepsilon (k)=1,\ \
\hat\varepsilon (e)=\hat\varepsilon (f)=0.
$$
The set $I$ of equivalence classes of irreducible
corepresentations of $(A,\D)$ is identified with the set $\2\7Z_+$
of non-negative half-integers. The fundamental corepresentation
($s=\2$) is defined by
\begin{equation}\label{e2.1}
U^\2 =(u^\2_{ij})_{ij}=
\pmatrix{\alpha & -q\gamma^* \cr
\gamma & \alpha^*},
\end{equation}
The corresponding representation of $\sl2$ is given by
\begin{equation} \label{e2.2}
k\mapsto\pmatrix{q^\2 & 0\cr 0 & q^{-\2}},\ \
e\mapsto\pmatrix{0 & 1 \cr 0 & 0}.
\end{equation}
The formulas (\ref{e2.1}--\ref{e2.2}) completely determine the
pairing between $(\A,\D)$ and $(\sl2,\hat\D)$. Note also that
\begin{equation}\label{e2.3}
\rho=f_{-1}=k^2.
\end{equation}

For $a\in\A$ and $\omega\in M(\hat\A)$ define elements in $\A$ by
$$
\pr_\omega(a)=\omega*a=(\iota\otimes\omega)\D(a)\ \ \hbox{and}\ \
a\r\omega=a*\omega=(\omega\otimes\iota)\D(a).
$$
Then the modular property of the Haar state $h$ can be expressed
as
\begin{equation}\label{e2.35}
h(a_1a_2)=h(a_2f_1*a_1*f_1)=h(a_2\pr_{k^{-2}}(a_1\r k^{-2})).
\end{equation}
For $n\in\7Z$ set
$$
\A_n=\{a\in\A\,|\,\pr_k(a)=q^{\frac{n}{2}}a\}.
$$
The norm closure $B=C(S^2_q)\subset A$ of $\A_0$ with the left
coaction $\D|_B$ of $(A,\D)$ is the quantum homogeneous sphere of
Podle\'{s}~\cite{P}. Note that $\B=\A_0$. The norm closure $A_n$
of $\A_n$ is an analogue of the space of continuous sections of
the line bundle over the sphere with winding number $n$.

\smallskip

It is well-known that the standard differential calculus on 
the quantum sphere~\cite{P1} can be obtained from various covariant 
differential calculi on $SU_q(2)$, see e.g.~\cite{S}.
Consider the $3D$-calculus of Woronowicz~\cite{Wor,KS}, that
is, the left-covariant first order differential calculus
$(\A\omega_{-1}\oplus\A\omega_0\oplus\A\omega_1,d)$ of $(\A,\D)$
with differential
$$
da=\pr_{fk^{-1}}(a)\omega_{-1}+
\frac{a-\pr_{k^{-4}}(a)}{1-q^{-2}}\omega_0+\pr_{ek^{-1}}(a)\omega_1
,
$$
and right $\A$-module action on the left-invariant forms
$\omega_{-1},\omega_0,\omega_1$ given by
$\omega_0a=\pr_{k^{-4}}(a)\omega_0$ and $\omega_i
a=\pr_{k^{-2}}(a)\omega_i$ for $i=\pm1$ and $a\in\A$. The
associated exterior differential algebra spanned by $a_0da_1\dots
da_n$, $a_i\in\A$, is completely described by the following rules:
$$
\omega_{-1}^2=\omega_0^2=\omega_1^2=0,\ \
\omega_{-1}\omega_1=-q^2\omega_1\omega_{-1},\ \
\omega_0\omega_1=-q^4\omega_1\omega_0,\ \
\omega_{-1}\omega_0=-q^4\omega_0\omega_{-1},
$$
$$
d\omega_{-1}=(q^2+q^4)\omega_0\omega_{-1},\ \
d\omega_0=-\omega_1\omega_{-1},\ \
d\omega_1=(q^2+q^4)\omega_1\omega_0.
$$
The standard differential calculus on $S^2_q$ is then obtained by
restricting the differential $d$ to $\B$, so the exterior algebra
$(\Gamma^\wedge_\B,d)$ for $S^2_q$ is the projective $\B$-module
given by $\Gamma^\wedge_\B=\span\{b_0db_1\cdots db_n\ |\
b_i\in\B\}$, and we have the following concrete description.

\begin{thm} \label{2.1}
Set $\Omega^{0,1}(\B)=\A_{-2}$, $\Omega^{1,0}(\B)=\A_2$ and
$\Omega^{1,1}(\B)=\B$. Then \enu{i} the de Rham complex on $S^2_q$
is the graded differential algebra
$$
\Gamma^\wedge_\B=\Gamma^{\wedge0}_\B\oplus\Gamma^{\wedge1}_\B
\oplus\Gamma^{\wedge2}_\B
=\B\oplus(\Omega^{0,1}(\B)\oplus\Omega^{1,0}(\B))
\oplus\Omega^{1,1}(\B)
$$
with multiplication $\wedge$ given by \newline$
(a_{0,0},\,a_{0,1},\,a_{1,0},\,a_{1,1})\wedge
(b_{0,0},\,b_{0,1},\,b_{1,0},\,b_{1,1})$
$$
=(a_{0,0}b_{0,0},\,a_{0,0}b_{0,1}+a_{0,1}b_{0,0},\,
a_{1,0}b_{0,0}+a_{0,0}b_{1,0},\,a_{0,0}b_{1,1}
-a_{0,1}b_{1,0}+q^2a_{1,0}b_{0,1}+a_{1,1}b_{0,0})
$$
and differential $d=\pr+\bar\pr$ given by
$$
\pr(a_{0,0},\,a_{0,1},\,a_{1,0},\,a_{1,1})
=(0,\,0,\,\pr_e(a_{0,0}),\,q\pr_e(a_{0,1})),
$$
$$
\bar\pr(a_{0,0},\,a_{0,1},\,a_{1,0},\,a_{1,1})
=(0,\,\pr_f(a_{0,0}),\,0,\,-q\pr_f(a_{1,0}));
$$
\enu{ii} if we define $\int\omega=h(\omega)$ for
$\omega\in\Gamma^{\wedge2}_\B$ by identifying
$\Gamma^{\wedge2}_\B$ with $\B$, then $\int$ is closed in the
sense that $\int d\omega=0$ for any
$\omega\in\Gamma^{\wedge1}_\B$; \enu{iii} the twist $\sigma$ of
$\B$ extends uniquely to an automorphism of $\Gamma^{\wedge}_\B$
commuting with~$d$, which we again denote by $\sigma$, and then
$$
\int\omega\wedge\omega'=(-1)^{\#\omega\#\omega'}
\int\sigma(\omega')\wedge\omega
$$
for any $\omega,\omega'\in\Gamma^\wedge_\B$ with
$\#\omega+\#\omega'=2$.
\end{thm}

\bp Part (i) follows from the equality
$\Gamma^\wedge_\B=\B+\A_{-2}\omega_{-1}+
\A_2\omega_1+\B\omega_1\omega_{-1}$, which is easily checked.
Since $\Gamma^{\wedge1}_\B=\B d\B$, to prove (ii) it is enough to
show that $\int db_1\wedge db_2=0$, that is
\begin{equation} \label{e2.4}
q^2h(\pr_e(b_1)\pr_f(b_2))=h(\pr_f(b_1)\pr_e(b_2)).
\end{equation}
We have
$$
h(\pr_e(b_1)\pr_f(b_2))=h(b_1\pr_{\hat
S(e)f}(b_2))=-q^{-1}h(b_1\pr_{ef}(b_2)).
$$
Similarly $h(\pr_f(b_1)\pr_e(b_2))=-qh(b_1\pr_{fe}(b_2))$. Since
$\pr_{ef}=\pr_{fe}$ on $\B$, assertion (ii) is proved. The
existence of an extension of $\sigma$ is obvious: it comes from
the automorphism $a\mapsto a\r k^2$ of $\A$. Concerning the
equality in part (iii), the only interesting case is when
$\omega\in\Omega^{1,0}(\B)$ and $\omega'\in\Omega^{0,1}(\B)$. In
this case we have to prove that for $a\in\A_2$ and $a'\in\A_{-2}$
we have
\begin{equation}\label{e2.5}
q^2h(aa')=h(\sigma(a')a).
\end{equation}
But this is true, since
$h(\sigma(a')a)=h(a\pr_{k^{-2}}(a'))=q^2h(aa')$ by
(\ref{e2.35}).\ep

Define the volume form by
$$
\tau(b_0,b_1,b_2)=\int b_0db_1\wedge
db_2=h(b_0(q^2\pr_e(b_1)\pr_f(b_2)-\pr_f(b_1)\pr_e(b_2))).
$$
Set also
$$
\tau_1(b_0,b_1,b_2)=h(b_0\pr_e(b_1)\pr_f(b_2))\ \ \hbox{and}\ \
\tau_2(b_0,b_1,b_2)=h(b_0\pr_f(b_1)\pr_e(b_2)).
$$
These forms were introduced in~\cite{SW}, where part of the
following proposition was proved.

\begin{prop} \label{2.2}
We have that\enu{i} $\tau$ is a cocycle in the twisted cyclic
complex $(C^\bullet_{\sigma,\lambda}(\B),b)$; \enu{ii} $\tau_1$
and $\tau_2$ are cocycles in the complex
$(\C^\bullet_{\sigma}(\B),b+B)$; \enu{iii} the cocycle
$q^2\tau_1+\tau_2\in\C^0_{\sigma}(\B)$ is a coboundary; in
particular, the cocycles $\tau$, $2q^2\tau_1$ and $-2\tau_2$ are
cohomologous.
\end{prop}

\bp Part (i) is a standard consequence of the properties of the
integral, cf~\cite{Co1,KMT}.

To prove (iii), set $\tilde\tau(b_0,b_1)=h(\pr_f(b_0)\pr_e(b_1))$.
Then
$$
(B_0\tilde\tau)(b_0)=(B^0_0\tilde\tau)(b_0)=\tilde\tau(b_0,1)
+\tilde\tau(1,\sigma(b_0))=0.
$$
On the other hand, using that $\pr_e$ and $\pr_f$ are
derivations on $\B$ and that $h(\sigma(b)a)=h(ab)$ for any
$a,b\in\B$, we get
$$
(b_2\tilde\tau)(b_0,b_1,b_2)
=h(b_0\pr_f(b_1)\pr_e(b_2))+h(\pr_f(\sigma(b_2))b_0\pr_e(b_1))
=\tau_2(b_0,b_1,b_2)+q^2\tau_1(b_0,b_1,b_2),
$$
where we have used (\ref{e2.5}) with $a=b_0\pr_e(b_1)$ and
$a'=\pr_f(a_2)$. Thus $q^2\tau_1+\tau_2$ is the coboundary
of~$\tilde\tau$. 

Since $\tau=\tau_1+\tau_2$ and $q^2\ne1$, statement (ii) follows
immediately from (i) and (iii). It can also be checked directly.
\ep

Note that although $\tau$ is cyclic, this is not the case for
$\tau_1$ and $\tau_2$.

\bigskip

\section{The Dirac Operator} \label{3}

The Dirac operator on the quantum sphere has been discovered,
studied and rediscovered in
various papers, see e.g. \cite{O,PS,DS,SW,M,Kr} and references 
therein.
We confine ourselves to summarize needed results in the 
next paragraph and in Proposition~\ref{3.1}.

One of the easiest ways to construct the Dirac operator is to
recall \cite{F} that the spinor bundle on a classical K\"{a}hler
manifold $M$ of complex dimension $m$ is
$\oplus^m_{i=0}\wedge^{0,i}\otimes{\cal S}$, where ${\cal S}$ is a
square root of the canonical line bundle $\wedge^{m,0}$, and then
the Dirac operator $D$ is $\sqrt{2}(\bar\pr+\bar\pr^*)$. In view
of the previous section these notions have rather
straightforward analogues for the quantum sphere. So consider the
space $H=L^2(A,h)$ of the GNS-representation of $A=C(SU_q(2))$
corresponding to the Haar state $h$. Then
$$
H=\oplus_{n\in\7Z}L^2(A_n,h).
$$
The left actions $a\mapsto\pr_\omega(a)$ and $a\mapsto a\r\hat
S^{-1}(\omega)$ of $\hat\A$ on $A$ extend to $*$-representations
of $\hat A$ on~$H$. These are, in fact, the left and the right
regular representations of $\hat A$. In the sequel we will
write~$\pr_\omega$ for the operators of the first representation
and simply $\omega$ for the operators of the second one. Then the
space of $L^2$-spinors and the Dirac operator are defined by
$$
H_+\oplus H_-=L^2(A_1,h)\oplus L^2(A_{-1},h)\ \ \hbox{and}\ \
D=\pmatrix{0 & \pr_e\cr \pr_f & 0},
$$
respectively. Note that $\pr_f(\A_n)\subset\A_{n-2}$ and
$\pr_e(\A_n)\subset\A_{n+2}$. This shows that $D$ is indeed an
operator on $H_+\oplus H_-$. It is worth noticing that when the
spinor bundle is given, the Dirac operator is essentially uniquely
determined by the requirements of $SU_q(2)$-invariance, boundedness
of commutators and the first order condition~\cite{DS}.

\begin{prop}\label{3.1}
We have \enu{i} $D^2=C$, where $\displaystyle
C=fe+\left(\frac{q^\2k-q^{-\2}k^{-1}}{q-q^{-1}}\right)^2$ is the
Casimir; \enu{ii} $\displaystyle [D,b]=\pmatrix{0 & q^\2\pr_e(b)
\cr q^{-\2}\pr_f(b) & 0}$ for any $b\in\B$.
\end{prop}

\bp We have $\displaystyle D^2=\pmatrix{\pr_{ef} & 0\cr 0 &
\pr_{fe}}$. Since $\pr_k|_{L^2(A_n,h)}=q^{\frac{n}{2}}$ and
$\pr_{ef}-\pr_{fe}=(q-q^{-1})^{-1}(\pr_{k^2}-\pr_{k^{-2}})$, we
see that $D^2=\pr_C$. To prove (i) it remains to note that
$\pr_C=C$ on $H$ as $C$ is in the center of $M(\hat\A)$ and $\hat
S(C)=C$.

Part (ii) follows from the identity
$$
[\pr_e,a]=\pr_e(a)\pr_{k^{-1}}+\pr_k(a)\pr_e-a\pr_e,
$$
which is valid for any $a\in\A$, and from a similar identity for $\pr_f$.
\ep

\medskip

We now want to develop a simple pseudo-differential calculus for
our Dirac operator. More precisely, we want to compute the
principal symbol of the operator of the form
$|D|^{2z}db|D|^{-2z}$, where $db=[D,b]$. For this it is more
convenient to consider a more general problem. The
Casimir~$C$ defines a scale of Hilbert spaces $\H_t$, $t\in\7R$,
so that for $t\ge0$
$$
\H_t=\hbox{Domain}(C^{t/2}),\ \ \|\xi\|_t=\|C^{t/2}\xi\|.
$$
Set $\H_\infty=\cap_t\H_t$. Note that $\H_\infty$ is a dense
subspace of $\H_t$ for every $t$. We say that an operator
$T\colon\H_\infty\to\H_\infty$ is of order $r\in\7R$ if it extends
by continuity to a bounded operator $\H_t\to\H_{t-r}$ for every
$t$. In this case we write $T\in op^r$. If $T\in op^0$, we denote
by $\|T\|_t$ the norm of the operator $T\colon\H_t\to\H_t$. Recall
that $\|T\|_t$ is a convex function of $t$. Note that $C^z\in
op^{2\Re z}$.

For $F\subset\7Z$, denote by $p_F$ the projection onto
$\oplus_{n\in F}L^2(A_n,h)$ with respect to the decomposition
$H=\oplus_{n\in\7Z}L^2(A_n,h)$.

\begin{prop} \label{3.2}
We have \enu{i} $\A\subset op^0$ and $p_F\in op^0$ for any
$F\subset\7Z$; \enu{ii} for $x=ap_F$ with $a\in\A_n$ and
$F\subset\7Z$ finite, the function $\7C\ni z\mapsto
x(z)=(C^zxC^{-z}-q^{nz}x)C^{\2}$ takes values in $op^0$ and has at
most linear growth on vertical strips; more precisely, for any
finite interval $I\subset\7R$, there exists $\lambda>0$ such that
$\|x(z)\|_t\le\lambda|z|$ for $t\in I$ and $\Re z\in I$.
\end{prop}

\bp Since $p_F$ is a projection commuting with $C$, it is obvious
that $p_F\in op^0$. It is also clear that to prove that $\A\subset
op^0$ it is enough to consider generators of $\A$. Similarly, if
(ii) is proved for $x_i=a_ip_{F_i}$, $a_i\in\A_{n_i}$, $i=1,2$,
and it is proved that $a_1,a_2\in op^0$, then
$x_1x_2=a_1a_2p_{F_1-n_2}p_{F_2}=a_1a_2p_{(F_1-n_2)\cap F_2}$ and
\newline
$\displaystyle (C^zx_1x_2C^{-z}-q^{(n_1+n_2)z}x_1x_2)C^{\2}$
$$
=(C^zx_1C^{-z}-q^{n_1z}x_1)C^{\2}(C^{-\2+z}x_2C^{\2-z})
+q^{n_1z}x_1(C^zx_2C^{-z}-q^{n_2z}x_2)C^{\2},
$$
so (ii) is also true for $a_1a_2p_{(F_1-n_2)\cap F_2}$.
Analogously, if (ii) is true for $ap_F$ with $a\in\A_n$, then it
is also true for $a^*p_{F+n}$. It follows that to prove (i) and
(ii) it is enough to consider generators of $\A$ and arbitrary
finite $F\subset\7Z$. We will only consider the generator
$\alpha$, the proof for $\gamma$ is similar. Consider the
orthonormal basis $\{\xi^s_{ij}\,|\,s\in\2\7Z_+,\
i,j=-s,\ldots,s\}$ given by normalized matrix coefficients,
$\xi^s_{ij}=q^id_s^{\2}u^s_{ij}\xi_h$, where $d_s=[2s+1]_q$. Then
$\alpha$ can be written as the sum of two operators $\alpha^+$ and
$\alpha^-$ such that
$$
\alpha^+\xi^s_{ij}=\alpha^{s+}_{ij}\xi^{s+\2}_{i-\2,j-\2}\ \
\hbox{and}\ \
\alpha^-\xi^s_{ij}=\alpha^{s-}_{ij}\xi^{s-\2}_{i-\2,j-\2}.
$$
What we have to know about the numbers $\alpha^{s+}_{ij}$ and
$\alpha^{s-}_{ij}$ is that they are of modulus $\le1$ and if
$|j|\le m$, then $|\alpha^{s+}_{ij}|\le\lambda_mq^s$ for some
$\lambda_m$ depending only on $m$ and $q$ (in fact, we have
$\alpha^{s+}_{ij}
=q^\2(d_sd_{s+\2})^{-\2}q^\frac{2s+i+j}{2}([s-i+1]_q[s-j+1]_q)^\2$,
cf~\cite{V}). In our basis we have also
$C\xi^s_{ij}=[s+\2]^2_q\xi^s_{ij}$. It follows immediately that
$\alpha^+,\alpha^-\in op^0$. Then

\medskip\noindent$\displaystyle
(C^{z+t}\alpha^-C^{-z-t}-q^zC^t\alpha^-C^{-t})C\xi^s_{ij}$
\begin{eqnarray*}
&=&[s+\2]^2_q\left(\left(\frac{[s]_q}{[s+\2]_q}\right)^{2z+2t}
  -q^z\left(\frac{[s]_q}{[s+\2]_q}\right)^{
  2t}\right)\alpha^{s-}_{ij}\xi^{s-\2}_{i-\2,j-\2}\\
&=&q^{z+t}[s+\2]^2_q\left(\frac{1-q^{2s}}{1-q^{2s+1}}\right)^{2t}
\left(\left(\frac{1-q^{2s}}{1-q^{2s+1}}\right)^{2z}-1\right)
  \alpha^{s-}_{ij}\xi^{s-\2}_{i-\2,j-\2}.
\end{eqnarray*}
Using the simple estimate $|(1-\beta)^z-1|=|\beta
z\int^1_0(1-\beta\tau)^{z-1}d\tau|\le\beta|z|(1-\beta_0)^{\Re
z-1}$ for $0\le\beta\le\beta_0<1$, we conclude that the operator
$(C^{z+t}\alpha^-C^{-z-t}-q^zC^t\alpha^-C^{-t})C$ is bounded, and
moreover, if $t$ and $\Re z$ lie in some finite interval then the
operator norm can be estimated by some multiple of $|z|$. Thus the
function $z\mapsto(C^{z}\alpha^-C^{-z}-q^z\alpha^-)C$ takes values
in $op^0$ and has at most linear growth on vertical strips. Then
surely the function
$z\mapsto(C^{z}\alpha^-C^{-z}-q^z\alpha^-)C^\2$ has the same
properties. Similarly the function
$z\mapsto(C^{z}\alpha^+C^{-z}-q^{-z}\alpha^+)C^\2$ takes values in
$op^0$ and has at most linear growth on vertical strips. Then

\medskip\noindent$\displaystyle
(C^{z}\alpha p_FC^{-z}-q^z\alpha p_F)C^\2$
$$
=(C^{z}\alpha^-C^{-z}-q^z\alpha^-)C^\2p_F
+(C^{z}\alpha^+C^{-z}-q^{-z}\alpha^+)C^\2p_F
+(q^{-z}-q^z)\alpha^+p_FC^\2.
$$
It remains to note that $p_F$ is the projection onto the space
spanned by the vectors $\xi^s_{ij}$ such that $-2j\in F$. Since
$|\alpha^{s+}_{ij}|\le\lambda_Fq^s$ if $-2j\in F$ for some
$\lambda_F$, we conclude that $\alpha^+p_FC^\2\in op^0$. \ep

It is worth noting that by analyzing the proof of Theorem~B.1
in~\cite{CM3} the slightly weaker result that $x(z)$ has
polynomial growth on vertical strip can be deduced using only the
fact that $Cx-q^{n}xC\in op^1$. This, in turn, follows from the
identity
$$
\pr_{fe}a-q^{n}a\pr_{fe}=\pr_{fe}(a)\pr_{k^{-2}}+\pr_{fk}(a)\pr_{k^{-1}e}
+\pr_{ke}(a)\pr_{fk^{-1}}
$$
for $a\in\A_n$. Thus a variant of Proposition~\ref{3.2} with
polynomial growth (which would, in fact, be sufficient for our
purposes) can be proved without any knowledge about the
Clebsch-Gordan coefficients.

\begin{cor}\label{3.3}
Consider the operator $\chi=\pmatrix{q & 0\cr 0 & q^{-1}}$ on
$H_+\oplus H_-$. Then for any $b\in\B$ there exists an analytic
function $z\mapsto b(z)\in B(H_+\oplus H_-)$ with at most linear
growth on vertical strips such that
$$
|D|^{-2z}db=db\chi^{2z}|D|^{-2z}+b(z)|D|^{-2z-1}
=\chi^{-2z}db|D|^{-2z}+b(z)|D|^{-2z-1}.
$$
\end{cor}

\bp Since $|D|^2=C$ and $\pr_f(b)\in\A_{-2}$, $\pr_e(b)\in\A_2$
for $b\in\B=\A_0$, this is an immediate consequence of
Propositions~\ref{3.1} and~\ref{3.2}. \ep

\medskip

Next we study the behavior of the heat operators. Consider the
direct sum of irreducible corepresentations of $(A,\D)$ with spins
$s\in\2+\7Z_+$. Denote by $\tilde H$ the space of this
corepresentation, and consider the corresponding representation of
$\hat A$ on $\tilde H$. Note that both $H_+=L^2(A_1,h)$ and
$H_-=L^2(A_{-1},h)$ can be identified with $\tilde H$ (since in
the notation of the proof of Proposition~\ref{3.2} the spaces
$H_+$ and $H_-$ are spanned by the vectors $\xi^s_{i,-\2}$ and
$\xi^s_{i,\2}$, respectively).

\begin{lem} \label{3.4}
On $\tilde H$ we have that \vspace{2mm}\enu{i} the function
$\eps\mapsto\eps\Tr(e^{-\eps C}\rho)$ is continuous and bounded on
$(0,\infty)$;\vspace{2mm}\enu{ii}
$\displaystyle\lim_{\eps\to0+}\eps(\Tr(e^{-\eps
C}\rho)-q^2\Tr(e^{-\eps q^2C}\rho))=0$; \enu{iii}
$\displaystyle\lim_{\eps\to0+}\eps\int^1_{q^2}\Tr(e^{-\eps
tC}\rho)dt=\mu=q^{-1}-q$.
\end{lem}

Using the Karamata theorem (see e.g. \cite{BGV}) it can be shown
that the limit $\displaystyle\lim_{\eps\to0+}\eps\Tr(e^{-\eps
C}\rho)$ does not exist.

\bpp{Lemma \ref{3.4}} If $f$ is a positive function then
$\Tr(f(C)\rho)=\sum^\infty_{n=1}[2n]_qf([n]^2_q)$. Since it is
obvious that the function $\eps\mapsto\eps\Tr(e^{-\eps C}\rho)$ is
continuous and vanishes at infinity, to prove (i) it is thus
enough to show that the function $\eps\mapsto\sum^\infty_{n=1}\eps
q^{-2n}e^{-\eps [n]^2_q}$ is bounded near zero. As
$[n]^2_q\ge(q^{-2n}-2)\mu^{-2}$, for $q^{2m}\le\eps\le q^{2(m-1)}$
we get
$$
\sum^\infty_{n=1}\eps q^{-2n}e^{-\eps [n]^2_q}
\le\sum^\infty_{n=1}q^{-2(n-m+1)}e^{-(q^{-2(n-m)}-2q^{2m})\mu^{-2}}
\le
q^{-2}e^{2\mu^{-2}}\sum^\infty_{n=-\infty}q^{-2n}e^{-q^{-2n}\mu^{-2}}.
$$
Thus (i) is proved. Since
$q^{-2n}\mu^{-2}\ge[n]^2_q\ge(q^{-2n}-2)\mu^{-2}$, we see also
that
$$
\lim_{\eps\to0+}\eps(\Tr(e^{-\eps C}\rho)-\sum^\infty_{n=1}
q^{-2n}\mu^{-1}e^{-\eps q^{-2n}\mu^{-2}})=0.
$$
Replacing $\eps$ by $\eps\mu^2$, to prove (ii) and (iii) we thus
have to show that
$$
\lim_{\eps\to0+}\eps(\sum^\infty_{n=1}q^{-2n}e^{-\eps
q^{-2n}}-q^2\sum^\infty_{n=1} q^{-2n}e^{-\eps q^2q^{-2n}})=0
$$
and
$$
\lim_{\eps\to0+}\eps\int^1_{q^2}dt\sum^\infty_{n=1}q^{-2n}e^{-\eps
tq^{-2n}}=1.
$$
Both statements are obvious. \ep

\bigskip

\section{Local Index Formula} \label{4}

Let us first briefly recall the proof of the local index formula
of Connes and Moscovici~\cite{CM3}, which is a far-reaching
generalization of the local index theorem of
Gilkey-Patodi~\cite{BGV}. Suppose we are given a $p$-summable even
spectral triple for an algebra $\B$. For each $\eps>0$ consider
the JLO-cocycle corresponding to the Dirac operator $\eps^\2D$, so
$$
\psi^\eps_{2n}(b_0,\ldots,b_{2n})=\eps^n\int_{\D_{2n}}dt
\Tr(\gamma b_0e^{-\eps t_0D^2}[D,b_1]e^{-\eps t_1
D^2}\ldots[D,b_{2n}]e^{-\eps t_{2n}D^2}),
$$
and note that we always have the following estimate
$$
|\psi^\eps_{2n}(b_0,\ldots,b_{2n})| \le
C_p\frac{\eps^{n-\frac{p}{2}}}{(2n)!}\Tr((1+D^2)^{-\frac{p}{2}})
\|b_0\|\prod^{2n}_{i=1}\|[D,b_i]\|,
$$
where $C_p$ is a universal constant. In the course of proving the
local index formula one provides a finite decomposition of
$\psi^\eps_{2n}$, $0\le 2n\le p$, of the form
\begin{equation} \label{e4.1}
\psi^\eps_{2n}=\sum_{k,l}\alpha^{(2n)}_{k,l}\eps^{-p_k(n)}(\log\eps)^l
+o(1),
\end{equation}
where $\alpha^{(2n)}_{k,l}$ are $2n$-cochains expressed in terms
of the residues of certain $\zeta$-functions. Since the pairing
with $K$-theory is independent of $\eps$, one concludes that it is
given by the pairing with the cocycle $(\alpha^{(2n)}_{0,0})_{0\le
2n\le p}$. The explicit form of this cocycle and the fact that it
is cohomologous to $(\Ch^{2n}(D))_n$ is the main content of the
local index formula.

To get (\ref{e4.1}) one needs to commute $e^{-\eps tD^2}$ through
$db$. For this one first uses the identity
$$
e^{-\eps tD^2}= \frac{1}{2\pi
i}\int_{C_\lambda}dz\Gamma(z)|D|^{-2z}(\eps t)^{-z},
$$
where $C_\lambda=\lambda+i\7R$, $\lambda>0$. Then one invokes the
expression
\begin{equation} \label{e4.2}
b_0|D|^{-2z_0}db_1|D|^{-2z_1}\ldots db_{2n}|D|^{-2z_n}
=\sum^{p-2n+1}_{k=0}f_k(b_0,\ldots,b_{2n};z_0,\ldots,z_{2n})
|D|^{-2(z_0+\ldots+z_{2n})-k}
\end{equation}
given by the pseudo-differential calculus for $D$, where $f_k$ are
analytic operator-valued functions with at most polynomial growth
on vertical strips. Since the last term ($k=p-2n+1$) is of
trace-class for $\Re z_i\ge\frac{2n-1}{2(2n+1)}$, it contributes
to $\psi^\eps_{2n}$ as $O(\eps^\2)$. On the other hand, assuming
that $\Tr(f_k(b_0,\ldots,b_{2n};z_0,\ldots,z_{2n})
|D|^{-2(z_0+\ldots+z_{2n})-k})$ extends to a nicely behaved
meromorphic function (note that a priori it is defined for $\Re
z_i\ge\frac{p-k}{2(2n+1)}$), the contribution of the $k$th term
can be estimated by counting the residues of this function in the
region $\frac{n}{2n+1}-\delta<\Re z_i<\frac{p-k}{2(2n+1)}$ for
some $\delta>0$.

\smallskip

We now want to apply this method to our situation. Although the
decomposition (\ref{e4.2}) will be different from the one obtained
in~\cite{CM3}, the method itself works perfectly well. So consider
the twisted entire cocycle $(\psi^\eps_{2n})_n$ given by
$$
\psi^\eps_{2n}(b_0,\ldots,b_{2n})
=\Ch^{2n}(\eps^{\2}D)(\rho;b_0,\ldots,b_{2n}).
$$
Using Lemma~\ref{1.2} we get
$$
\psi^\eps_0(b_0)=\Tr(\gamma b_0e^{-\eps
D^2}\rho)=h(b_0)\Tr(e^{-\eps C}\rho)-h(b_0)\Tr(e^{-\eps C}\rho)=0,
$$
where we think of the Casimir $C$ as the operator acting on
$\tilde H$, see the end of Section~\ref{3}.
Consider~$\psi^\eps_2$. Fix $b_0,b_1,b_2\in\B$ and set
$$
\zeta(z_0,z_1,z_2)=\Tr(\gamma
b_0|D|^{-2z_0}db_1|D|^{-2z_1}db_2|D|^{-2z_2}\rho).
$$
This function is well-defined and analytic in the region
$\Re(z_0+z_1+z_2)>1$. By Corollary~\ref{3.3},
$$
\zeta(b_0,b_1,b_2)=\Tr(\gamma
b_0db_1db_2\chi^{2z_1}|D|^{-2(z_0+z_1+z_2)}\rho)+
\tilde\zeta(z_0,z_1,z_2),
$$
where $\displaystyle \tilde\zeta(z_0,z_1,z_2)=$
$$
=\Tr(\gamma
b_0db_1\chi^{2z_0}b_2(z_0+z_1)|D|^{-2(z_0+z_1+z_2)-1}\rho) +
\Tr(\gamma b_0b_1(z_0)|D|^{-2(z_0+z_1)-1}db_2|D|^{-2z_2}\rho).
$$
The function $\tilde\zeta$ is holomorphic in the region
$\Re(z_0+z_1+z_2)>\2$. Moreover, if we fix
$\lambda,\lambda_0>\frac{1}{6}$, $\lambda>\lambda_0$, then
$$
\tilde\zeta(z_0,z_1,z_2)(1+|z_0|)^{-1}(1+|z_1|)^{-1}
$$
is bounded in the region $\lambda_0<\Re z_i<\lambda$. Fix
$\lambda>\frac{1}{3}$ and $\frac{1}{6}<\lambda_0<\frac{1}{3}$.
Then
$$
\psi^\eps_2(b_0,b_1,b_2)=\frac{\eps}{(2\pi
i)^3}\int_{\D_2}dt\int_{C_\lambda^3}dz
\Gamma(z_0)\Gamma(z_1)\Gamma(z_2)(\eps t_0)^{-z_0} (\eps
t_1)^{-z_1}(\eps t_2)^{-z_2}\zeta(z_0,z_1,z_2),
$$
where $C_\lambda=\lambda+i\7R$. If we consider the same expression
for $\tilde\zeta$, we can replace integration over $C_\lambda^3$
by integration over $C_{\lambda_0}^3$, since the integrand is
holomorphic and vanishes at infinity. It follows that such an
expression is $O(\eps^{1-3\lambda_0})$. Thus, as $\eps\to0$,
\begin{eqnarray*}
\psi^\eps_2(b_0,b_1,b_2)&=&\frac{\eps}{(2\pi
i)^3}\int_{\D_2}dt\int_{C_\lambda^3}dz
\Gamma(z_0)\Gamma(z_1)\Gamma(z_2)(\eps t_0)^{-z_0} (\eps
t_1)^{-z_1}(\eps t_2)^{-z_2}\\
& &\hspace{1cm}\times\Tr(\gamma
b_0db_1db_2\chi^{2z_1}|D|^{-2(z_0+z_1+z_2)}\rho)+o(1)\\
&=&\eps\int_{\D_2}dt\Tr(\gamma
b_0db_1db_2e^{-\eps(t_0+\chi^{-2}t_1+t_2)D^2}\rho)+o(1)\\
&=&\eps\int_{\D_2}dt\Tr(\gamma
b_0db_1db_2e^{-\eps(1+(\chi^{-2}-1)t_1)D^2}\rho)+o(1)\\
&=&\eps\int^1_0(1-t)\Tr(\gamma
b_0db_1db_2e^{-\eps(1+(\chi^{-2}-1)t)D^2}\rho)dt+o(1).
\end{eqnarray*}
By Proposition~\ref{3.1}(ii), we have
$$
b_0db_1db_2=\pmatrix{b_0\pr_e(b_1)\pr_f(b_2) & 0\cr 0 &
b_0\pr_f(b_1)\pr_e(b_2) }.
$$
Hence by Lemma~\ref{1.2} we get
\begin{eqnarray}
\psi^\eps_2(b_0,b_1,b_2)&=&h(b_0\pr_e(b_1)\pr_f(b_2))\eps\int^1_0(1-t)
\Tr(e^{-\eps(1+(q^{-2}-1)t)C}\rho)dt\nonumber\\
& &\hspace{1cm}-h(b_0\pr_f(b_1)\pr_e(b_2))\eps\int^1_0(1-t)
\Tr(e^{-\eps(1+(q^2-1)t)C}\rho)dt+o(1). \label{e4.10}
\end{eqnarray}
We have
\begin{equation} \label{e4.11}
\eps\int^1_0(1-t)\Tr(e^{-\eps(1+(q^2-1)t)C}\rho)dt
=-\frac{\eps}{(1-q^2)^2}\int^1_{q^2}(q^2-t)\Tr(e^{-\eps
tC}\rho)dt.
\end{equation}
On the other hand,
\begin{eqnarray}
\eps\int^1_0(1-t) \Tr(e^{-\eps(1+(q^{-2}-1)t)C}\rho)dt &=&
\frac{\eps}{(1-q^2)^2}\int^1_{q^2}(1-t)\Tr(e^{-\eps
q^{-2}tC}\rho)dt\nonumber\\
&=&\frac{\eps q^2}{(1-q^2)^2}\int^1_{q^2}(1-t)\Tr(e^{-\eps
tC}\rho)dt+o(1) \label{e4.12}
\end{eqnarray}
by Lemma \ref{3.4}(ii). By Lemma \ref{3.4}(iii),
$\eps\int^1_{q^2}\Tr(e^{-\eps tC}\rho)dt\to q^{-1}(1-q^2)$. Thus
putting (\ref{e4.10}-\ref{e4.12}) together we get
$$
\psi^\eps_2=\frac{q}{1-q^2}(\tau_1+\tau_2)
-(q^2\tau_1+\tau_2)\frac{\eps}{(1-q^2)^2}\int^1_{q^2}\Tr(e^{-\eps
tC}\rho)tdt+o(1).
$$
Finally, consider $\psi^\eps_{2n}$ with $n>1$. Combining
(\ref{e1.10}) with Lemmas 10.8 and 10.11 in \cite{GBVF} yields the
following standard estimate
$$
|\psi^\eps_{2n}(b_0,\ldots,b_{2n})| \le
C_p\frac{\eps^{n-\frac{p}{2}}}{(2n)!}\Tr(C^{-\frac{p}{2}}\rho)
\|b_0\|\prod^{2n}_{i=1}\max_{0\le t\le 1}\{\|[D,b_i\r\rho^t]\|\},
$$
where $p$ is any number larger than $2$. Thus we get the first
part of the following theorem, which is the main result of the
paper.

\begin{thm} \label{4.1}
For $\eps>0$, consider the cocycle $(\psi^\eps_{2n})_n$ in the
twisted entire cyclic cohomology of $\B$ given by
$$
\psi^\eps_{2n}(b_0,\ldots,b_{2n})
=\Ch^{2n}(\eps^{\2}D)(\rho;b_0,\ldots,b_{2n}).
$$
Then
\begin{eqnarray*}
\psi^\eps_0&=&0;\\
\psi^\eps_2&=&\frac{q}{1-q^2}(\tau_1+\tau_2)
-(q^2\tau_1+\tau_2)\frac{\eps}{(1-q^2)^2}\int^1_{q^2}\Tr(e^{-\eps
tC}\rho)tdt+o(1)\ \ \hbox{as}\ \ \eps\to0;\\
\|\psi^\eps_{2n}\|&\le& C_\delta\frac{\eps^{n-1-\delta}}{(2n)!}\ \
\hbox{for any}\ \ \delta>0\ \ \hbox{and}\ \ n>1,
\end{eqnarray*}
where the norm of a multi-linear form on $\B$ is defined using the
norm
$$
\max_{0\le t\le 1}\{\|b\r\rho^t\|+\|[D,b\r\rho^t]\|\}
$$
on $\B$. In particular, the map $q{\rm-Ind}_D$ on
$K^{SU_q(2)}_0(C(S^2_q))$ is given by the pairing with the twisted
cyclic cocycle $-q^{-1}\tau$, where
$\tau(b_0,b_1,b_2)=h(b_0(q^2\pr_e(b_1)\pr_f(b_2)
-\pr_f(b_1)\pr_e(b_2)))$.
\end{thm}

\bp The pairing of $(\psi^\eps_{2n})_n$ with $K$-theory does not
depend on $\eps$. Using the Karamata theorem it can be shown that
the limit $\lim_{\eps\to0+}\eps\int^1_{q^2}\Tr(e^{-\eps
tC}\rho)tdt$ does not exist. It follows that the cocycle
$q^2\tau_1+\tau_2$ pairs trivially with $K$-theory. Thanks to
Proposition~\ref{2.2}(iii) we know even more, $q^2\tau_1+\tau_2$
is a coboundary. We also conclude that the pairing is given by the
cocycle $\frac{q}{1-q^2}(\tau_1+\tau_2)\in\C^0_\sigma(\B)$, which
is cohomologous to the cocycle $\frac{1}{2q}\tau$. Recalling
definition (\ref{e1.3}) of the map
$C^2_{\sigma,\lambda}(\B)\to\C^0_\sigma(\B)$, we see that the
pairing is defined by the twisted cyclic cocycle $-q^{-1}\tau$.
\ep

We end the paper with an actual computation of indices. Observe
first that the spaces $A_n$ can be considered as equivariant
Hilbert $B$-modules, and thus they define elements of
$K_0^{\hat\A}(B)$ which we denote by $[A_n]$. As in the classical
case, the group $K_0^{\hat\A}(B)$ is a free abelian group with
generators $[A_n]$, see~\cite{NT}. Consider now the module $A_1$
and note that $\A_1=\alpha\B+\gamma\B$. The map $T\colon
H_\2\otimes B\to A_1$ given by
$$
T(\xi_{-\2}\otimes b)=q\gamma b\ \ \hbox{and}\ \ T(\xi_\2\otimes
b)=-\alpha b,
$$
where $H_\2$ is the space of the spin $\2$ corepresentation of
$(A,\D)$, see Section~\ref{2}, is an equivariant partial isometry.
Thus $A_1$ is isomorphic to the equivariant Hilbert $B$-module
$p(H_\2\otimes B)$ with projection
$$
p=T^*T=\pmatrix{q^2\gamma^*\gamma & -\alpha\gamma^*\cr
-\gamma\alpha^* & \alpha^*\alpha}.
$$
The explicit form of a projection corresponding to $A_n$ for
arbitrary $n$ can be found in~\cite{HM}. Let $\Tr_s$ be the trace
on $\hat\A$ defined by the spin $s$ representation of $\hat\A$,
and let $\phi_s$ be the corresponding $q$-trace,
$\phi_s=\Tr_s(\cdot\rho)$. By definition we have, for any twisted
cyclic cocycle $\varphi\in C^2_\sigma(\B)$, that
$$
\<[\varphi],[p]\>=\sum_{i_0,i_1,i_2\atop
j_0,j_1,j_2}\phi_\2(m_{i_0j_0}m_{i_1j_1}m_{i_2j_2})
\varphi(p_{i_0j_0},p_{i_1j_1},p_{i_2j_2})
=\sum_{i_0,i_1,i_2}q^{-2i_0}\varphi(p_{i_0i_1},p_{i_1i_2},p_{i_2i_0}).
$$
Using the formula
$h((\gamma^*\gamma)^n)=(1-q^2)(1-q^{2(n+1)})^{-1}$ for the Haar
state, a lengthy but straightforward computation yields
$$
q{\rm-Ind}_D([A_1])=\indq([A_1])=\<[-q^{-1}\tau],[p]\>=-1.
$$
This is enough to conclude that the equivariant index
$\indq([A_1])$ equals $-\Tr_0$. To see this we shall use a
continuity argument for $q\in(0,1)$. Write $\alpha(q)$,
$\gamma(q)$, and so on, to distinguish operators for different
$q$. The spaces $L^2(C(SU_q(2)),h)$ can be identified for all $q$.
We also identify the spaces $H_+\oplus H_-$ of $L^2$-spinors. Note
that $F=D(q)|D(q)|^{-1}$ is independent of $q$ (in the notation of
the proof of Proposition~\ref{3.2} we have
$F\xi^s_{i,-\2}=\xi^s_{i,\2}$). The functions $q\mapsto\alpha(q)$
and $q\mapsto\gamma(q)$ are norm-continuous as can easily be
verified by looking at the Clebsch-Gordan coefficients. It follows
that our quantum Bott projections $p(q)\in B(H_\2\oplus H_+\oplus
H_-)$ depend continuously on~$q$. Let $I_s(q)\in\hat\A(q)$ be the
support of the spin $s$ representation. Considered as operators on
$H_\2\oplus H_+\oplus H_-$ the projections $I_s(q)$ depend
continuously on $q$ (and are, in fact, finite-rank operators). As
the functions
$$
m_s(q)=(2s+1)^{-1}{\rm Ind}(p(q)_-I_s(q)(1\otimes F)p(q)_+I_s(q))
$$
are continuous and integer-valued, they are constant. We have by
definition
$$
\ind([p(q)])=\sum_sm_s\Tr_s.
$$
Since
$$
-1=\ind([p(q)])(\rho)=\sum_sm_s[2s+1]_q,
$$
and the functions $q\mapsto[n]_q$, $n\in\7N$, are linearly
independent on any infinite set, we conclude that $m_0=-1$ and
$m_s=0$ for $s>0$. Thus $\ind([p])=-\Tr_0=-\hat\eps$. This, in
turn, is sufficient in order to find the non-equivariant Chern
character.

\begin{prop}
The image of the non-equivariant Chern character of our Fredholm
module in $HP^0(\B)$ coincides with the class of the cyclic
$0$-cocycle $\tau'$ given by
$$
\tau'(\alpha^{n-m}\gamma^m{\gamma^*}^n) =\cases{(1-q^{2n})^{-1}\ \
\hbox{for}\ \ n=m>0,\cr 0\hspace{19mm}\hbox{otherwise},}
$$
where we used the convention $\alpha^k=(\alpha^*)^{-k}$ for $k<0$.
In particular, $\ind([A_n])(1)=-n$.
\end{prop}

\bp The cocycle $\tau'$ was found in \cite{MNW2}, and is one of
the two generators of $HP^0(\B)\cong\7C^2$. Since the class of a
cocycle in $HP^0(\B)$ is completely determined by its pairing with
$[1]$ and $[p]$, we conclude that the Chern character is
cohomologous to $\tau'$. The equality $\<[\tau'],[A_n]\>=-n$ was
established in~\cite{H}.\ep

The fact that the non-equivariant Chern character is cohomologous
to a $0$-cocycle is natural as our spectral triple is
$\eps$-summable for any $\eps>0$. On the other hand, the spectral
triple is $(2+\eps,\rho)$-summable in the sense of~\cite{NT}, so
twisted cyclic cohomology does not see the dimension drop and
captures the volume form.

\smallskip

We finally remark that $\ind([A_n])=-{\rm
sign}(n)\Tr_\frac{|n|-1}{2}$ for $n\ne0$. To prove this it
suffices to check that $\<[-q^{-1}\tau],[A_n]\>=-[n]_q$. Another
possibility is to use the classical theory. To this end one just
has to show that there are projections $p_n(q)$ representing
$[A_n]$ with the property that $I_s(q)p_n(q)$ depend continuously
on $q\in(0,1]$ (the projections $p_n(q)$ themselves can be
discontinuous at $q=1$). In the classical case the operator
$p_n(1\otimes D)p_n$ is homotopic to the operator $D_n=\pmatrix{0
& \pr_e\cr \pr_f & 0}$ which acts on the Hilbert space
$L^2(A_{n+1},h)\oplus L^2(A_{n-1},h)$. Both operators are
differential operators of order~$1$ with the same principal
symbol, and the index of $\pr_f\colon L^2(A_{n+1},h)\to
L^2(A_{n-1},h)$ is given by the Borel-Weil-Bott theorem and can
also easily be found by direct computations.

\bigskip

\flushleft {Sergey Neshveyev, Mathematics Institute, University of
Oslo, PB 1053 Blindern, Oslo 0316, Norway\\
{\it e-mail}: neshveyev@hotmail.com}

\flushleft {Lars Tuset, Mathematics Institute, University of
Oslo, PB 1053 Blindern, Oslo 0316, Norway\\
{\it e-mail}: Lars.Tuset@iu.hio.no}

\end{document}